\documentclass{amsart}

\usepackage{amssymb}
\usepackage{amsmath}
\usepackage{array}

\usepackage[small,nohug,heads=vee]{diagrams}
\diagramstyle[lablestyle=\scriptstyle]

\newcommand{\Pf}{{\em Proof}. }
\newcommand{\EPf}{\hfill$\square$}

\newcommand\liegr{\sf}

\newcommand{\SU}[1]{\mbox{${\liegr SU}(#1)$}}

\newcommand{\U}[1]{\mbox{${\liegr U}(#1)$}}
\newcommand{\SP}[1]{\mbox{${\liegr Sp}(#1)$}}
\newcommand{\SO}[1]{\mbox{${\liegr SO}(#1)$}}

\newcommand{\OG}[1]{\mbox{${\liegr O}(#1)$}}
\newcommand{\Spin}[1]{\mbox{${\liegr Spin}(#1)$}}
\newcommand{\G}{\mbox{${\liegr G}_2$}}
\newcommand{\F}{\mbox{${\liegr F}_4$}}
\newcommand{\E}[1]{\mbox{${\liegr E}_{#1}$}}

\newcommand\fieldsetc{\mathbb}

\newcommand{\Z}{\fieldsetc{Z}}
\newcommand{\R}{\fieldsetc{R}}
\newcommand{\C}{\fieldsetc{C}}
\newcommand{\Q}{\fieldsetc{H}}

\newtheorem*{thm}{Theorem}

\newtheorem{lem}{Lemma}

\title{Representations with $Sp(1)^k$-reductions and 
quaternion-K\"ahler symmetric spaces}

\author[C.~Gorodski]{Claudio Gorodski}
\author[F.~J.~Gozzi]{Francisco J.~Gozzi}

\address{Instituto de Matem\'atica e Estat\'\i stica, Universidade de
S\~ao Paulo, Rua do Mat\~ao, 1010, S\~ao Paulo, SP 05508-090, Brazil}

\email{gorodski@ime.usp.br}

\thanks{The 
first author has been partially supported 
by the CNPq grant 303038/2013-6 and the FAPESP project 2011/21362-2.}

\address{Instituto de Matem\'atica e Estat\'\i stica, Universidade de
S\~ao Paulo, Rua do Mat\~ao, 1010, S\~ao Paulo, SP 05508-090, Brazil}

\email{fj.gozzi@gmail.com}

\thanks{The second author has been supported by
  the FAPESP fellowship 2014/22568-1.}

\date{\today}

\subjclass[2010]{57S15, 22E46, 53C35}

\begin{document}

\begin{abstract}
We classify non-polar irreducible representations of connected compact Lie
groups whose orbit space is 
isometric to that of a
representation of a finite extension of $\SP1^k$
for some $k>0$. It follows that they are obtained from
isotropy representations of certain quaternion-K\"ahler symmetric
spaces by restricting to the ``non-$\SP1$-factor''.
\end{abstract}

\maketitle


\section{Introduction}

The aim of this paper is to contribute to the program initiated
in~\cite{GL}, namely, hierarchize the  
representations of compact Lie groups in terms of the
complexity of their orbit spaces, viewed as metric spaces.
We say that two representations are \emph{quotient-equivalent}
if they have isometric orbit spaces. Given a representation, if there is a 
quotient-equivalent representation of a lower-dimensional group, then we say
that the former representation \emph{reduces} to the latter one and that
the latter representation is a \emph{reduction} of the former one.  
A \emph{minimal reduction} of a representation is a 
reduction with smallest
possible dimension of the underlying group. 

At the basis of the hierarchy lie 
the \emph{polar representations}, namely, those representations
that reduce to finite group actions, and which turn out to be 
related to symmetric spaces (indeed every polar representation of a
connected compact Lie group has the same orbits as the isotropy representation
of a symmetric space~\cite{D}). 
In~\cite{GL2}, there were studied
and classified those irreducible representations of connected groups
that reduce to actions of groups
with identity component a torus $S^1\times\cdots\times S^1$. It was 
shown that, mostly, those are close relatives of Hermitian symmetric spaces. 
Herein we study the quaternionic version, namely, those irreducible 
representations of connected groups
that reduce to an action of a group whose identity component is a 
``quaternionic torus'' $S^3\times\cdots\times S^3$. Interestingly enough,
these are related to quaternion-K\"ahler symmetric spaces, in a stricter
sense than in the Hermitian case.  

\begin{thm}
  Let $\tau: H \to \OG W $ be a non-polar
  irreducible representation of a connected
compact Lie group. Assume that $\tau$ is quotient-equivalent
to a representation $\rho:G\to \OG V$ where
$G^0\cong\SP1^k$ for some $k>0$ and $\dim G<\dim H$. 
Then $k=3$, $G$ is disconnected and $V=\otimes^3\C^2$; moreover,
the cohomogeneity of $\tau$ is $7$ and it is obtained by restricting 
the isotropy representation of a certain quaternion-K\"ahler symmetric space
to the ``non-$\SP1$-factor''. More precisely, $\tau$ is one of:
\[ \begin{array}{|c|c|c|}
\hline
\tau & \mbox{qK symmetric space} & G/G^0\\
\hline
(\SO n\times \SP1,\R^n\otimes_{\mathbb R}\Q) & \SO{n+4}/(\SO n\times\SO4) & \sf A_1\\
(\SP3,\Lambda^3\C^6\ominus\C^6) & \F/(\SP3\SP1) & \sf A_2\\
(\SU6,\Lambda^3\C^6) & \E6/(\SU6\SU2) & \sf A_2\\
(\Spin{12},\C^{32}) & \E7/(\Spin{12}\SU2) & \sf A_2\\
(\E7,\C^{56}) & \E8/(\E7\SU2) & \sf A_2\\
\hline
\end{array}\]
Finally, $\rho$ is a minimal reduction of~$\tau$. 
\end{thm}

Note that there are three additional families of quaternion-K\"ahler
symmetric spaces absent from the Table, namely, the 
restrictions of their isotropy representations do not have the kind of reduction
as in the Theorem. The Theorem says in particular that $k$ must be 
odd, so $\rho$ is a representation of quaternionic type. 
A posteriori we find that the $\SP1$-group of isometries
of $W/H=V/G$ induced by the
normalizer of $G$ in $\OG V$ can be lifted to a $\SP1$-subgroup
of $\OG W$. In this connection, we note that
the general problem of lifting isometries
of orbit spaces to isometries deserves further
attention (\cite[Question 1.6]{AL} and~\cite[Question 1.13]{GL}).

The structure of the proof of the Theorem goes as follows.
The situation of a non-trivial reduction from~$\tau$ to $\rho$ entails
the presence of non-empty boundary for the orbit space of $\rho$,
according to results in~\cite[\S5]{GL}. A careful analysis
of the existence of boundary points then forces $k=3$. Now
$\rho|_{G^0}$ is of quaternionic type, so its orbit space admits
an $\SP1$-group of isometries. An application of a theorem
of Thorbergsson, as in~\cite{GL2}, shows that $\tau$
must be the restriction of a polar representation of Coxeter type
$\sf B_4$ or $\sf F_4$. The argument is finished by invoking
Dynkin's classification of maximal connected closed subgroups of
compact Lie groups.

We follow notation and terminology from~\cite{GL}. 
See also~\cite{GL3} for some background material on stratification
of orbit spaces and their metric structures. The authors wish to thank
Alexander Lytchak for several valuable comments. 

\section{Structure of the examples}

In this section we show that the representations $\tau$
listed in the table of the Theorem indeed
admit reductions to a group whose identity component
is $\SP1^3$. In all cases, the reduction in hand is the
Luna-Richardson-Straume reduction to the normalizer
of a principal isotropy group acting on the fixed point set
of this principal isotropy group~\cite[\S2.6]{GL3}. 
It will follow from the discussion
in Subsection~\ref{conn} that those reductions are minimal.

\subsection{The real Grassmannian} The isotropy representation
of the rank $4$ real Grassmannian manifold is
$\SO n\times\SO 4$ acting on $\R^n\otimes_{\mathbb R}\R^4$.
In this case $\tau$ is given by $H=\SO n\times\SU2$ acting
on $\R^4\otimes_{\mathbb R}\C^2$ with principal isotropy
group $H_{princ}=\SO{n-4}$, and its effective normalizer
$N_H(H_{princ})/H_{princ}=\mathsf{S}(\OG4\times\OG{n-4})\times\SU2/H_{princ}=\OG4\times\SU2=\mathbb Z_2\cdot \SP1^3$ has the desired form.

\subsection{The exceptional cases} For each one of the rank $4$
exceptional quaternion-K\"ahler symmetric spaces,
the isotropy representation $\hat\tau$ is given by
$\hat H=H\cdot\SP1$ acting on $\hat W=W\otimes_{\mathbb H}\Q\cong W$ and
$\tau=\hat\tau|_H:H\to\OG W$ is a representation of quaternionic type. 
The principal isotropy
group of $\hat\tau$ has the form $\hat H_{princ}=H_{princ}\cdot Q$
where $Q$ is the quaternion group $\{\pm1,\pm i,\pm j,\pm k\}$
diagonally embedded in $H\cdot \SP1$. Since $\hat\tau$ is
asystatic (see~\cite[pp.11-12]{str} or~\cite[\S2.2]{GK}), it is polar and
the fixed point set $W^{\hat H_{princ}}$ is a section;
this is a totally real subspace of $W$ of dimension $4$.
Now $W^{H_{princ}}$ is a quaternionic subspace of $W$, and it must
be the quaternionic span of $W$, of real dimension $16$.
Since the cohomogeneity of $\tau$ is 
$c(\tau)=c(\hat\tau)+\dim\SP1=4+3=7$ (compare \cite[Table~A]{HH}), 
we deduce that 
$\dim N_H(H_{princ})/H_{princ}=\dim W^{H_{princ}}-c(\tau)=16-7=9$.
The group $\dim N_H(H_{princ})/H_{princ}$ acts irreducibly on
$W^{H_{princ}}$, so its center is at most one-dimensional.
A quick enumeration of the possible groups reveals that
$\dim N_{H}(H_{princ})/H_{princ}$ is locally isomorphic to
$\SP1^3$ or $\U3$, but $W^{H_{princ}}$ is a representation
of quaternionic type, and the latter group admits none.
Note also that the only $16$-dimensional irreducible representation of
$\SP1^3$ is $\otimes^3\C^2$. 

\section{Proof of the Theorem}

Throughout this section we let $\tau: H \to \OG W $ and $\rho:G\to \OG V$
be as in the statement of the Theorem.

\subsection{Triviality of the principal isotropy group}

We observe that the principal isotropy group of $\rho^0:=\rho|_{G^0}$ 
is trivial. Indeed we have the following general result.

\begin{lem}
Every irreducible representation of $\SP1^k$ is either polar
or has trivial principal isotropy group. 
\end{lem}

\Pf We shall collect the irreducible
representations with possibly non-trivial principal isotropy
group and check that they are all polar. 
The complexification of the principal isotropy group of a representation 
of $\SP1^k$ on $V$ is the stabilizer in general position (sgp)
of the corresponding complex representation of $\mathsf{SL}(2,\mathbb C)^k$
on~$V^c=V\otimes_{\mathbb R}\C$ \cite[\S5]{S}. If $V$ has no invariant
complex structure, then $V^c$ is irreducible; otherwise, $V^c$ consists
of two copies of $V$. Now it suffices to start listing irreducible
representations $U$ of $\mathsf{SL}(2,\C)$ such that:
\begin{itemize}
\item $U$ has non-trivial sgp, if $U$ is a representation of real type;
\item $U\oplus U$ has non-trivial sgp, if $U$ is a representation of
  quaternionic type;
\end{itemize}
and then checking that these admit a invariant real form~$V$ which is a
polar representation of~$\SP1^k$. 

At this point, one could apply~\cite[Theorem]{A-E-V} and refer
to \cite[Table~1]{P1} and \cite[Table~1]{P2} to finish the job.
Alternatively, the case $k=1$ is well known (see e.g.~\cite[p.~231]{PV};
here $V=\R^3$ or $\R^5$ which are polar).
For $k\geq2$, we work from scratch using~\cite[Theorem~2]{AP} (see
also~\cite[Theorem~7.10]{PV})
to deduce that a necessary condition for non-triviality of sgp of $V^c$ is that
$\dim V=\dim_{\mathbb C} V^c\leq8k$; this estimate takes care of almost
all representations. Additional representations are excluded with the help
of the Corollary to Lemma~6 in~\cite{AP}. We end up with few
non-polar irreducible
representations of $\SP1^k$, each of which is easily
checked to have trivial principal
isotropy group. \EPf

\subsection{General setting}\label{gen}

Note that $V$ equals $\C^{n_1}\otimes_{\mathbb C}\cdots\otimes\C^{n_k}$ 
or a real form $[\C^{n_1}\otimes_{\mathbb C}\cdots\otimes\C^{n_k}]_{\mathbb R}$ 
according to whether the dimension $n_i$ is even for an odd, resp.~even, 
number of indices~$i$. We have that $V$ is of quaternionic type in the 
first case, and of real type in the second one. 
We may assume $2\leq n_1\leq\cdots\leq n_k$. 

Due to~\cite[\S5]{GL},
$X=V/G=W/H$ has non-empty boundary as an Alexandrov space. 
We recall that the boundary is given by the closure of the codimension 
one strata of $X$. Since the principal isotropy is trivial, 
a point $p\in V$ projects to such a stratum --- called a $G$-\emph{important
point} in~\cite{GL} ---
if and only if its isotropy group $G_p$ is a sphere $S^\ell$ of dimension~$\ell$
equal to~$0$, $1$ or~$3$ and 
\begin{equation}\label{dim-form}
 \dim V - 1 - \ell = 3k-m+f,
\end{equation}
where $m$ is the dimension of the normalizer $N_G(G_p)$ 
and $f$ is the dimension of the fixed point set $V^{G_p}$ of $G_p$ in $V$, 
see ~\cite[Lemma~4.1]{GL}.

In the following two subsections, we will analyze the boundary
of $V/G$ and show that it can be non-empty only if $G$ is disconnected,
$k=3$ and $V=\otimes^3\C^2$. The cases $k=1$ and $k=2$ were discussed
in~\cite[\S10]{GL},
so we may assume $k\geq3$.

\subsection{Connected case}\label{conn}
We first suppose $G$ is connected and show that this assumption is
incompatible with $V/G$ having non-empty
boundary.   
So suppose that $p\in V$ is a $G$-important point. 
Since $G$ is connected, no $G$-important point 
may lie in an exceptional orbit~\cite{Lyt}, so $G_p$ 
is not discrete. 
Moreover, any $\SU2$-subgroup of $G$ contains a unique involution that 
is central in $G$. Since such an involution cannot have fixed points
by irreducibility, we cannot have $G_p\cong\SU2$. We deduce that 
$G_p\cong S^1$. 

The dimension formula~(\ref{dim-form}) yields
\[ \theta\cdot n_1\cdots n_k - 2 =3k-m+f, \]
where $\theta=1$ or $2$
in case $V$ is of real or quaternionic type, respectively.
It is moreover clear that $m\geq k$. Since $V^{G_p}$, resp. 
its complexification, is the sum of weight spaces whose weights
lie in a hyperplane of the dual Lie algebra of the maximal torus
of $G$, it is not hard to see that $f\leq \theta\cdot n_2\cdots n_k$. 
We deduce that 
\[ \theta(n_1-1)n_2\cdots n_k\leq 2k+2. \]
In particular, $\theta\cdot 2^{k-2}\leq k+1$ implying $k=3$ or~$4$. 
If $k=3$, we obtain $(\SP1^3,\otimes^3\C^2)$ and
$(\SO4\times\SO3,\R^4\otimes_{\mathbb R}\R^3)$; the latter
representation is polar, which cannot be.  
If $k=4$, we obtain 
$(\SO4\times\SO4,\R^4\otimes_{\mathbb R}\R^4)$ which is also polar
and again is excluded. 

It remains to analyze the case of $(\SP1^3,\otimes^3\C^2)$. 
Here the dimension formula~(\ref{dim-form}) says 
that $5+m=f$, and $m\geq3$ implying $f\geq8$. Recall that the weights
of $\C^2$ are $\pm\epsilon$, where $2\epsilon$ is the positive root
of $\SP1$. It is apparent that $V^{G_p}$ can contain at most four weight
spaces so $f=8$. However, in this case $G_p$ is a circle diagonally embedded
in two factors of $\SP1^3$ which gives $m=5$ and contradicts the 
dimension formula. We deduce that $G$ cannot be connected, as 
desired.  

We have shown that the orbit space of any non-polar representation of $\SP1^k$
has empty boundary. Together with~\cite[Proposition~5.2]{GL}, this implies
that \emph{an irreducible action of any extension of $\SP1^k$ by
a finite group is either polar or reduced (i.e.~it cannot be further reduced)}.

\subsection{Disconnected case}\label{disconn}
We now suppose $G$ is disconnected and prove that 
$V/G$ can have non-empty boundary
only if $k=3$ and $V=\otimes^3\C^2$.
Since $H$ is connected, 
there is an involution $w\in G\setminus G^0$, called a 
\emph{nice involution}, that acts as a reflection
on $V/G^0$ (see Proposition~3.2 and~\S4.3 in~\cite{GL}).
The dimension formula reads
\begin{equation}\label{dim-form2}
 \dim V - 1 =\dim G - \dim C(w) + \dim V^w 
\end{equation}
where $C(w)$ is the centralizer of $w$ in $G$ and $V^w$ 
is the fixed point set of $w$ in $V$.
The element~$w$ acts on $G^0$ by conjugation.  

\subsubsection{Outer automorphism}\label{outer} We first discuss the case 
in which $w$ acts on $G^0$ by an outer
automorphism. Since $w$ is not inner, its action on $G^0$ 
induces a non-trivial involutive permutation~$\sigma$ of the factors,
and its action on~$V$ induces a corresponding permutation of the factors of $V$.
Consider $w_0\in\OG V$ given by
\[ w_0(v_1\otimes\cdots\otimes v_k)=v_{\sigma(1)}\otimes\cdots\otimes v_{\sigma(k)}. \]
Then $w_0^{-1}w$ induces an inner automorphism of $G^0$, so we can write
$w=w_0hz$ where $h=(h_1,\ldots,h_k)\in G^0$ and $z\in\OG V$ centralizes $G^0$. 

Consider $V$ as a representation~$\rho^0$ of $G^0$ and relabel the indices
to write
\[ V=V_1\otimes_{\mathbb R}\cdots\otimes_{\mathbb R} V_a\otimes_{\mathbb R}\R^{n_1}\otimes_{\mathbb R}
\cdots\otimes_{\mathbb R}\R^{n_b}\otimes_{\mathbb R}\Q^{p_1/2}\otimes_{\mathbb H}
\cdots\otimes_{\mathbb F} \Q^{p_c/2} \]
where $n_i\geq3$ is odd, $p_j\geq2$ is even, $\mathbb F=\Q$ or $\R$ 
according to whether $c$ is even or odd;
moreover,  $V_i=\R^{m_i}\otimes_{\mathbb R}\R^{m_i}$ or 
$V_i=\Q^{m_i/2}\otimes_{\mathbb R}\Q^{m_i/2}$ according to whether
$m_i$ is odd or even, and $w_0$ exchanges the factors of $V_i$ for $i=1,\ldots,a$,
and fixes the other factors of $V$.
If~$c$ is even, then $V$ is of real type, so $z=\pm1$; if $c$ 
is odd, then $V$ is of quaternionic type, so $z$ is right multiplication
on $\Q^{p_c/2}$ by some element of~$\SP1$
(we view quaternionic vector spaces as right $\Q$-modules); 
in the latter case, $w_0$ fixes the factor $\Q^{p_c/2}$, so in any case
$w_0$ commutes with $z$.
Now $w^2=1$ gives that
\[ (h_{\sigma(1)}h_1,\ldots,h_{\sigma(k)}h_k)=z^{-2}\in Z(G^0)=\{\pm1\}^k. \]
We deduce that $h_{2i}=\pm h_{2i-1}^{-1}$ for $i=1,\ldots,a$.
Now we can take
\[ \tilde h=(1,h_1^{-1},\ldots,1,h_\ell,1,\ldots,1)\in G^0, \]
where $\ell=2a-1$, to replace $w$ by the conjugate element
$\tilde w=\tilde h w \tilde h^{-1}$ in somehow simpler form, namely,
\[ \tilde w=w_0z(1,\pm1,\ldots,1,\pm1,h_{\ell+1},\ldots,h_k) \]
(we could also have some simplification for the $h_i$ for $i>\ell$, but it
is unimportant to the sequel). 

We have
\[ \dim V =  m_1^2\cdots m_a^2 \cdot n_1\cdots n_b\cdot p_1\cdots p_c \cdot 
\theta, \]
where $\theta=1$ or~$2$ whether $c$ is even or odd, 
\[ \dim G = 3(2a+b+c)\quad\mbox{and}\quad \dim C(w)=\dim C(\tilde w)\geq3a+b+c. \]
The $\pm1$-eigenspaces of~$w_0$ on $V_i$ have dimension 
$\frac{m_i(m_i\pm1)}2=:m_i^{\pm}$.
We deduce that 
\[ \dim V^w =\dim V^{\tilde w}\leq \dim V^{\pm w_0}\leq\dim V^{w_0} = M\cdot
n_1\cdots n_b\cdot p_1\cdots p_c\cdot\theta \]
where
\[ M^{even}=\sum m_1^{\pm}\cdots m_a^{\pm} \] 
and the sum runs through all combinations with an even number of negative
signs. 
The dimension formula~(\ref{dim-form2}) now gives
\[ m_1\cdots m_an_1\cdots n_bp_1\cdots p_c\theta\left(
m_1\cdots m_a-M^{even}\right)
\leq 3a+2b+2c+1. \]
The factor between parenthesis is $M^{odd}=\sum m_1^{\pm}\cdots m_a^{\pm}$ where
the combinations are now taken with an odd number of negative signs.
Since $m_i^+\geq\frac32$ and $m_i^-\geq\frac12$, we estimate
\begin{eqnarray*}
  M^{odd} & = & \sum_{\ell=0}^{\left[\frac{a-1}2\right]}\binom{a}{2\ell+1}
  \frac{3^{a-(2\ell+1)}}{2^a} \\
  &=& \left(\frac32\right)^a\frac12\left[\left(1+\frac13\right)^a-
    \left(1-\frac13\right)^a\right] \\
  &=&\frac12\left(2^a-1\right)\\
  &\geq& 2^{a-2}.
\end{eqnarray*}
We deduce that 
\[ 2^{2a-2}3^b2^{c+1}\cdot\theta\leq 3a+2b+2c+1. \]
It immediately follows that $a\leq2$, 
and a run through the possibilities,
excluding polar representations, yields that 
$\rho^0$ must be $(\SP1^3,\otimes^3\C^2)$.

\subsubsection{Inner automorphism}
We next consider the case in which $w$ acts on~$G^0$ as an inner 
automorphism and show that this case gives nothing. 
Write  $w=qj$ where $q$ centralizes $G^0$ and $j\in G^0$.
Write also
\[  V=\R^{m_1}\otimes_{\mathbb R}
\cdots\otimes_{\mathbb R}\R^{m_a}\otimes_{\mathbb R}\Q^{n_1/2}\otimes_{\mathbb H}
\cdots\otimes_{\mathbb F} \Q^{n_b/2} \]
where $m_i\geq3$ is odd and $n_j\geq2$ is even. Suppose $V$ is of real 
type ($b$ is even). Then $q=\pm1$. Since $q$ does not lie in $G^0$, we must 
have $q=-1$ and $b=0$, namely, all factors of $G^0$ are 
isomorphic to $\SO3$ and $w=-j$, where $j\in G^0$, $j^2=1$.
Write $j=j_1\cdots j_a$ where 
$j_i$ is the component of~$j$ in the $i$th factor of $G^0$, and 
assume $j_i\neq1$ precisely for $1\leq i\leq a'$ for some $0\leq a'\leq a$.   
On one hand, the dimension formula~(\ref{dim-form2}) gives
\begin{eqnarray*}
 \dim V^j&=&\dim V - \dim V^w\\
         &=& \dim G - \dim C(w) + 1 \\
 &=& 3a-(a'+3(a-a'))+1 \\
 &=& 2a'+ 1.
\end{eqnarray*}
On the other hand, 
\begin{eqnarray*}
\dim V^j & = & m_{a'+1}\cdots m_{a}\dim V^{j_{1}\cdots j_{a'}} \\
   &=& m_{a'+1}\cdots m_{a} \sum g_{1}\cdots g_{a'} 
\end{eqnarray*}
where $g_i=e_i$ or~$f_i$ 
and the sum runs through all possibilities
with an even number of $f_i$'s; here
\[ e_i = \dim V^{j_i}=\left\{\begin{array}{ll}
p_i&\mbox{if $p_i$ is odd,}\\
p_i+1&\mbox{if $p_i+1$ is even}
\end{array}\right. \]
and 
\[ f_i = \dim V^{-j_i}=\left\{\begin{array}{ll}
p_i&\mbox{if $p_i+1$ is odd,}\\
p_i+1&\mbox{if $p_i$ is even}
\end{array}\right. \]
and $m_i=2p_i+1\geq3$ for all $i$. 

Since $e_i\geq1$ and $f_i\geq2$, we estimate
\begin{eqnarray*}
\sum g_{1}\cdots g_{a'} &=& e_{1}\cdots e_{a'}+ f_{1}f_{2}e_{3}\cdots e_{a'}+\cdots\\
&\geq& \sum_{\ell=0}^{[\frac{a'}2]}\binom{a'}{2\ell}2^{2\ell}\\
&=&\frac12[3^{a'}+(-1)^{a'}].
\end{eqnarray*}
We deduce that 
\[ 2a'+1\geq m_{a'+1}\cdots m_a\frac12[3^{a'}+(-1)^{a'}]\quad\mbox{and} 
\quad 2a'+1\geq \frac12[3^{a'}+(-1)^{a'}]. \]
The second inequality can be satisfied only if $a'\leq2$;
in this case, using $a\geq3$, we see that the first inequality is never
satisfied. Hence $V$ cannot be of real type. 

We finally take up the case $V$ is of quaternionic type ($b$ is odd). 
Then $q^2=j^{-2}$ lies in the center of $G^0$, which is 
isomorphic to $(\mathbb Z_2)^b$, and $q$ is a unit quaternion multiplying on the right. 
It follows that $q^2=\pm1$. If $q^2=1$, then irreducibility of $G^0$ implies
that one of the $\pm1$-eigenspaces of $q$ is trivial, namely, $q=\pm1$, but then 
$q\in G^0$ as $b$ is odd, a contradiction. If $q^2=-1$, then $q$ defines a 
complex structure with respect to which $w$ is a complex involution.
In particular, the fixed point set $V^w$ has even dimension.  
Note that $m:=\dim C(w)\in\{a+b,a+b+2,\ldots,a+3b\}$ implying $m\equiv a+b\ \mathrm{mod}\ 2$. 
The dimension formula yields
\[ \dim V - 1 = 3(a+b) - m + \dim V^w \] 
leading to a contradiction as $\dim V$ is even, too. \EPf 

\subsection{Discrete part}
It follows from the discussion in the 
previous two subsections that $G\neq G^0=\SP1^3$
and $V=\otimes^3\C^2$. The next step is 
to determine the exact nature of $G$.

\begin{lem}\label{discrete}
The group $G=\Gamma\ltimes G^0$ 
where $\Gamma=\sf A_1$ or $\sf A_2$.
\end{lem}

\Pf Recall that the Coxeter group $\sf A_n$ is the
group of permutations on $n+1$ elements.  
We saw in Subsection~\ref{disconn} that 
there exists a nice involution $w\in G\setminus G^0$
acting on $G^0$ by an
outer automorphism. Note that $\dim V=16$, $\dim G=9$ and $m=6$ 
since $\sigma$ permutes two factors of $G^0$. The dimension formula
implies that $f=12$ and hence $w$ operates
on $\otimes^3\C^2$ by interchanging two factors.
Since $G$ acts with trivial principal isotropy groups, 
$G/G^0\to\mathrm{Isom}(V/G^0)$ 
is injective; moreover, its image is a group generated by reflections. 
We select a  generating set consisting of nice involutions $w\in G$. 
Since they all have the form above, it follows that they generate a 
subgroup $\Gamma$ of $G$, and thus $G=\Gamma\ltimes G^0$, where $\Gamma$ is 
a subgroup of $\sf A_2$. Note that 
$\Gamma$ can only be $\sf A_1$ or $\sf A_2$, 
since it is non-trivial and contains 
an element of order $2$. \EPf

\medskip

Both cases described in Lemma~\ref{discrete} do occur: 
the first one for the rank $4$
real Grassmannian and the second one for the rank $4$ exceptional 
quaternion-K\"ahler 
symmetric spaces. We will see below that no further examples exist.

\subsection{Symmetries of $X$}
The key point is that $\rho^0:\SP1^3\to\OG V$ is a representation
of quaternionic type. The centralizer $\SP1'$ of $G^0$ in $\OG V$ 
also centralizes $\Gamma$ and acts thus 
on $X=V/G=W/H$:
\begin{diagram}
&W&&&\\
&\dTo&&&\\ 
X=&W/H&=\otimes^3\C^2/\Gamma\cdot\SP1^3&\curvearrowleft\SP1'\\
&&\dTo&\\
&&Y=\otimes^2\R^4/\Gamma\cdot\SO4^2\\
\end{diagram}
Denote the composite map $W\to Y$ by~$\pi$. Since 
$(\SO4\times\SO4,\R^4\otimes_{\mathbb R}\R^4)$ is a polar representation
(indeed the isotropy representation of the rank~$4$ real Grassmannian
manifold),
$Y$ is a flat Riemannian orbifold (of dimension~$4$) and hence 
the components of the level sets of~$\pi$  
yield an isoparametric foliation $\mathcal F$ by full irreducible 
submanifolds (compare~\cite[\S2]{GL2}). The codimension of $\mathcal F$ 
is $4$, so by a theorem of Thorbergsson~\cite{Th2}, $\mathcal F$ is homogeneous, 
namely, the maximal connected subgroup $\hat H$ of $\OG V$ that preserves the
leaves of $\mathcal F$ acts transitively on them. By definition, $\hat H$ is closed,
acts polarly, and contains~$H$. In particular, it acts irreducibly on $W$. It follows 
that $\hat\tau:\hat H\to \OG W$ is
the isotropy representation of an
irreducible symmetric space, of rank $4$. 

The geometry of $Y$ can be understood. The Coxeter group of
the polar representation 
$(\SO4\times\SO4,\R^4\otimes_{\mathbb R}\R^4)$ is $\sf D_4$, so
$Y=\R^4/\Gamma'$ where $\Gamma'$ is a finite extension
of $\sf D_4$ by $\Gamma$, 
of order $2\#(\sf D_4)$ or $3!\#(\sf D_4)$,
that acts irreducibly on $\R^4$. 
Recall that $\mathrm{Aut}(\sf D_4)\subset\mathrm{Ad}(\sf F_4)$,
that is, every automorphism
of $\sf D_4$ becomes inner in~$\sf F_4$ (see e.g.~\cite[Theorem~14.2]{Ad}).
Since the representation of $\Gamma'$ on $\R^4$ is irreducible of real type,
we deduce that every element of $\Gamma'$ differs from an element of $\sf F_4$
by $\pm 1$, but $-1\in\sf D_4$. This proves that $\Gamma'$ is a subgroup of
$\sf F_4$. Since $\sf F_4=\sf D_4\ltimes\sf A_2$ and 
$\sf B_4=\sf D_4\ltimes\sf A_1$, we deduce 
that $\Gamma'=\sf B_4$ or $\sf F_4$ according to
whether $\Gamma=\sf A_1$
or $\sf A_2$.

\subsection{End of the proof}

We refer to the classification of isotropy representations
of symmetric spaces~\cite[ch.~8]{wolf}
to list the possibilities for $\hat\tau$. If $\Gamma'=\sf B_4$, then
$\hat\tau:\hat H\to \OG{W}$ is one of:
\[ \begin{array}{cc}
(\SO n\times\SO4,\R^n\otimes_{\mathbb R}\R^4)\ \quad\mbox{($n\geq5$)}, & (\U8,\Lambda^2\C^8), \\
   ({\sf S}(\U n\times\U4),\C^n\otimes_{\mathbb C}\C^4)\ \quad\mbox{($n\geq4$)}, & (\U9,\Lambda^2\C^9), \\
  (\SP n\times\SP4,\Q^n\otimes_{\mathbb H}\Q^4)\ \quad\mbox{($n\geq4$)},  & (\U4,\rm S^2\C^4). \end{array} \]
If $\Gamma'=\sf F_4$, then
$\hat\tau$ is associated to one of the exceptional 
quaternion-K\"ahler symmetric spaces of rank $4$, namely, the last four
cases in the table of the Theorem. In each case, we look
for closed subgroups $H$ of $\hat H$ that act irreducibly with
cohomogeneity~$7$ on~$W$. A straightforward, not very long calculation
using the lists of maximal connected
closed subgroups of compact Lie groups compiled by Dynkin~(\cite{Dyn}; see
also~\cite{gp})
reduces the possibilities for $H$ to those listed in the Theorem
plus two extra candidates, namely,  
$\Spin7\times\U2\subset\SO8\times\SO4$ and
$\SU4\times\SP2\subset\mathsf{S}(\U4\times\U4)$, both corresponding to the case
$\Gamma'=\sf B_4$. In the sequel we rule out those representations
by showing that they cannot be quotient-equivalent to
$(\Z_2\cdot\SP1^3,\otimes^3\C²)=(\OG4\times\SP1,\R^4\otimes_{\mathbb R}\R^4)$.

Owing to the discussion in Subsections~\ref{conn} and~\ref{disconn},
the boundary of the orbit space $X$ of $(\OG4\times\SP1,\R^4\otimes_{\mathbb R}\R^3)$
originates from the fixed point set of a nice involution
$w\in(\OG4\setminus\SO4)\times\{1\}$, which we can take to
be~$w=\mathrm{diag}(-1,1,1,1)$. It follows that $\partial X$ is also given as the
orbit space of a representation (and we say that $X$ has
\emph{linear boundary}), namely,
$(\OG3\times\SP1,\R^3\otimes_{\mathbb R}\R^4)$. It has already been
remarked in Subsection~\ref{gen} that this representation has empty boundary.
We deduce that $X$ contains no strata of codimension~$2$ along its boundary. 
Another remark that will be useful below is that
the slice representations at non-zero points of 
$(\OG3\times\SP1,\R^3\otimes_{\mathbb R}\R^4)$ 
are all infinitesimally polar.
Recall that a representation of a compact Lie group is called
\emph{infinitesimally polar} if the slice representations
at non-zero points are all polar~\cite{GL3};
in the case at hand, this property follows because the isotropy subgroups
of $\OG3\times\SP1$ can be at most $1$-dimensional.
Note that infinitesimal polarity for a representation can be detected
metrically as it is equivalent to having the orbit space isometric to a
Riemannian orbifold~\cite{LT}. 

It follows easily from~\cite[\S13]{S} that the orbit space of
$(\SU4\times\SP2,\C^4\otimes_{\mathbb C}\Q^2)$ admits no $S^3$-boundary components
(and it admits no $\Z_2$-boundary components, as the underlying group is connected),
but we can find two $S^1$-boundary components that meet at a codimension $2$
strata, namely, the projections of the fixed point sets of the circle
subgroups
\[ (\mathrm{diag}(e^{i\theta},e^{i\theta},e^{-i\theta},e^{-i\theta}),(
e^{-i\theta},e^{-i\theta},e^{i\theta},e^{i\theta})) \]
and
\[ (\mathrm{diag}(e^{i\theta},1,e^{-i\theta},1),(
e^{-i\theta},1,e^{i\theta},1)), \]
as is readily checked. Hence this orbit space cannot be isometric
to~$X$.

Finally, we turn to $(\Spin7\times\U2,\R^8\otimes_{\mathbb R}\C^2)$
and its orbit space, which we denote by $Z$. 
The slice representation at a point $p$
given by a pure tensor, say restricted to
the identity component, is given
by~$(\G\times\SO2,\R^7\otimes_{\mathbb R}\R^2\oplus\R^7)$ which, in view of the classification in~\cite{GL3}, is not infinitesimally polar. On the other hand, the point~$p$ projects to a $S^3$-boundary
component of~$Z$ given by the projection of the fixed point set of a 
$\SU2$-subgroup of $\SU3\subset\G\times\{1\}\subset\G\times\SO2$.
Hence~$Z$ cannot be isometric to $X$.

That $\rho$ is a minimal reduction of~$\tau$ follows from the last
assertion in Subsection~\ref{conn} and this finishes
the proof of the Theorem. 

\providecommand{\bysame}{\leavevmode\hbox to3em{\hrulefill}\thinspace}
\providecommand{\MR}{\relax\ifhmode\unskip\space\fi MR }
\providecommand{\MRhref}[2]{%
  \href{http://www.ams.org/mathscinet-getitem?mr=#1}{#2}
}
\providecommand{\href}[2]{#2}


\end{document}